\documentclass[a4paper, 12pt]{article}
\usepackage{amssymb,amsbsy,amsmath,amsfonts,amssymb,amscd}
\usepackage{latexsym,eucal,exscale, epsfig}

\newcommand{\pk}{{\mathfrak p}_K}
\newcommand{\ofr}{\mathfrak o}
\newcommand{\pfr}{\mathfrak p}
\newcommand{\bk}{{\mathbf k}}
\newcommand{\ds}{\displaystyle}

\newcommand{\tX}{{\tilde X}}
\newcommand{\lra}{\longrightarrow}
\newcommand{\KK}{{\mathcal K}}
\newcommand{\HH}{\mathcal H}
\newcommand{\BB}{\mathbf B}
\newcommand{\TT}{\mathbf T}
\newcommand{\GG}{\mathbf G}
\newcommand{\UU}{\mathbf U}
\newcommand{\hh}{\mathfrak h}

\newcommand{\VV}{{\mathcal V}}
\newcommand{\noi}{\noindent}
\newcommand\cind{{\hbox{\rm $c$-ind}}}
\newcommand\St{{\mathbf S}{\mathbf t}}
\newcommand\tdX{\tilde X}
\newcommand\tdY{\tilde Y}
\newcommand\CC{\mathbb C}
\newcommand\ZZ{\mathbb Z}
\newcommand{\matrice}[4]{\left(\begin{array}{cc} #1 & #2 \\
#3 & #4 \end{array}\right)}
\title{A geometric construction of types for  the smooth
  representations of ${\rm PGL}(2)$ of a local field}
\author{Paul Broussous}

\date{March 2012}

\begin{document}
\maketitle

\abstract{We show that almost all (Bushnell and Kutzko) types of ${\rm PGL}(2,F)$, $F$
 a non-Archimedean locally compact field of odd residue characteristic, naturally 
appear in the cohomology of finite graphs.}

\setcounter{section}{0}

\vskip2cm

\centerline{\bf Introduction} 
\bigskip

Let $F$ be a non-Archimedean locally compact field and $G$ be the
group ${\rm PGL}(2,F)$. We assume that {\it the residue characteristic of
$F$ is not $2$}. In previous works (\cite{[Br1]}, \cite{[Br2]})
we defined a tower of directed graphs $(\tdX_n )_{n\geqslant 0}$ lying
$G$-equivariantly over the Bruhat-Tits tree $X$ of $G$. We proved the
two following facts :
\bigskip

\noi {\bf Theorem 1} (\cite{[Br2]}, Theorem (3.2.4), page 502). {\it
  Let $(\pi ,\VV )$ be a non-spherical generic smooth irreducible
  representation. Then $(\pi ,\VV )$ is  a quotient of the cohomology
  space with compact support $H_{c}^{1}(\tdX_{n(\pi )}, \CC )$, where
  $n(\pi )$ is the conductor of $\pi$.}
\bigskip

\noi {\bf Theorem 2} (\cite{[Br2]}, Theorem (5.3.2), page 512). {\it
  If $(\pi ,\VV )$ is supercuspidal smooth irreducible representation
  of $G$, then we have :
$$
{\rm dim}_{\CC}\, {\rm Hom}_G\, [H_{c}^{1}(\tdX_{n(\pi )}, \CC
  ), \VV ] = 1\ .
$$}
\bigskip

In this paper we make the $G$-module structure of $H_c^{1}(\tdX_n ,\CC
)$ more explicit  for all $n\geqslant 0$, and draw some
interesting consequences.
\bigskip

 Let  us fix an edge $[s_0 ,s_1 ]$ of $X$ and denote by $\KK_0$ and
 $\KK_1$ the stabilizers in $G$ of $s_0$ and $[s_0,s_1 ]$
 respectively. Then $\KK_0$ and $\KK_1$ form a set of representatives of
the two  conjugacy classes of maximal compact subgroups in $G$. If $n$
is even, we have a $G$-equivariant mapping $p_n$~:
$\tdX_n \lra X$ which respects the graph structures. 
We denote by $\Sigma_n$ the subgraph $p_n^{-1}([s_0
  ,s_1 ])$. If $n$ is odd, then after passing to the first barycentric
subdivisions, we have a $G$-equivariant  mapping $p_n$~:
$\tdX_n \lra X$ which respects the graph structures. 
We denote by $\Sigma_n$ the subgraph $p_n^{-1} (S(s_0
,1/2))$, where $S(s_0 ,1/2 )$ denotes the set of points $x$ in $X$
such that $d(x,s_0 )\leqslant 1/2$ (here $d$ is the natural distance
on the standard geometric realization of 
$X$, normalized in such a way that $d(s_0 ,s_1 )=1$). 

Then for all $n$, $\Sigma_n$ is a finite graph, equipped with a
an action of $\KK_1$ if $n$ is even, and $\KK_0$ if $n$ is odd. So the
cohomology spaces $H^1 (\Sigma_n ,\CC )$ provide finite dimensional
smooth representations of $\KK_1$ or $\KK_0$, according to the parity
of $n$.

 For each $n\geqslant 0$, we define an finite set ${\mathcal P}_n$
 of pairs $(\KK ,\lambda )$ formed of a maximal compact subgroup
 $\KK\in \{ \KK_0 ,\KK_1\}$ and of an irreducible smooth
 representation of $\KK$. By definition we have $(\KK ,\lambda )\in
 {\mathcal P}_n$ if and only if there exists $k\in \{ 0,1,...,n\}$
 such that $(\KK ,\lambda )$ is an irreducible constituent of the
 representation    $H^1 (\Sigma_k ,\CC )$. For $(\KK ,\lambda )\in {\mathcal P}_n$ 
and $k\leqslant n$, we denote by $m_\lambda^k$ the multiplicity of 
$\lambda$ in $H_c^1 (\Sigma_k ,\CC )$ and we set $m_{n,\lambda}= m_\lambda =n_\lambda^0 + \cdots n_\lambda^n$. 
Note that $n_\lambda$ depends on $(\KK ,\lambda )$ and $n$.

\bigskip

 The main results of this article are the following.
\bigskip

\noi {\bf Theorem A}. {\it For all $n\geqslant 0$, we have the direct sum
decomposition :
$$
H_{c}^{1}(\tdX_n ,\CC )= \St_G\oplus \bigoplus_{(\KK ,\lambda )\in
  {\mathcal P}_n} (\cind_{\KK}^{G}\ \lambda )^{m_\lambda} \ .
$$}

\noi (Here $\St_G$ denotes the Steinberg representation of $G$).
\bigskip

\noi {\bf Theorem B}. {\it For all $n\geqslant 0$, any element of
  ${\mathcal P}_n$ is} 
\bigskip

\noi a) {\it either a type in the sense of Bushnell and Kutzko's
  type theory \cite{[BK]}, which is not a type for the unramified
  principal series}
\smallskip

\noi b) {\it or a pair of the form $(\KK_0 ,\chi\circ {\rm det}\otimes
  \St_{\KK_0})$, where $\chi$ is a smooth character of $F^{\times}$ of
  order $2$, trivial on the group of $1$-units in $F^\times$, 
 and $\St_{\KK_0}$ is the representation inflated from the
  Steinberg representation of ${\rm PGL}(2)$ of the residue field of $F$},
\smallskip

\noi c) {\it or the pair $(\KK_1 , {\mathbf 1}_{\KK_1} )$, 
where $\mathbf 1$ denotes a trivial character.}
\bigskip

\noi {\bf Corollary C}. {\it Let $n\geqslant 0$. If 
$(\KK ,\lambda )\in {\mathcal P}_n$ is a cuspidal type, 
then  $m_{\lambda ,n} =1$. }
\bigskip

 Indeed this follows from Theorems 2 and A using Frobenius 
reciprocity for compact induction.
\bigskip

 By Theorem 1, any Bernstein component of $G$, different from the unramified 
principal series component, must have a type in ${\mathcal P}_n$ 
for $n$ large enough.   {\it Hence the grapĥs $\tdX_n$, $n\geqslant 0$,  
provide a geometric construction of types for almost all Bernstein components of $G$.}
\smallskip

 We conjecture that if $(\KK ,\lambda )\in {\mathcal P}_n$ is a type of $G$, 
then  $n_\lambda =1$. 
\medskip

 Finally let us observe that this construction gives a new proof that
 the irreducible supercuspidal
 representations of $G$ are obtained by compact induction. Our proof
  differs from  original Kutzko's proof (\cite{[Ku]},
 also see \cite{[BH]}) only at the exhaustion steps. Indeed our
 ``supercuspidal'' types
 are the same as Kutzko's, but we prove that any irreducible
 supercuspidal representation contains such a type by using an
 argument based on \cite{[Br1]} and \cite{[Br2]}, that is mainly on
 the existence of the new vector.
\bigskip

The article is organized as follows. The proof of Theorem A relies first on
combinatorial properties of the graphs $\tdX_n$ that are stated and proved
in {\S}2. Using this combinatorial properties and some homological arguments, 
we show in {\S}3 how to relate the cohomology  of $\tdX_n$ to that of $\tdX_{n-1}$.
 The irreducible components of $H^1 (\Sigma_n )$ are determined in {\S}4 when 
$n$ is even, and in {\S}5 and {\S}6 when $n$ is odd. A synthesis of the
arguments of paragraphs 2 to 6, leading to a proof a theorem A and B, is given in {\S}7.

\bigskip

We shall assume that the reader is familiar with the language of
Bushnell and Kutzko's type theory \cite{[BKtypes]} and with the
language of strata (\cite{[BK]}, \cite{[BH]}).

\section{Notation}

We shall denote by 
\medskip

$F$ a non-Archimedean non-discrete locally compact field,

$\ofr$ its valuation ring,

$\pfr$ the maximal ideal of $\ofr$,

$\varpi$ the choice of a uniformizer of $\ofr$,

$\bk=\ofr /\pfr$ the residue field of $F$,

$p$ the characteristic of $\bk$,

$q=p^f$ the cardinal of $\bk$,

$G$ the group ${\rm PGL}(2,F)$.

$t_{\varpi}$ the image of the matrix $\matrice{1}{0}{0}{\varpi}$ in
$G$. 
\bigskip

 The results of this article are obtained under the
\bigskip

\noi {\bf Hypothesis}. {\it The characteristic of $\bk$ is not $2$}
\bigskip

We shall often define an element,  a subset, or a subgroup of $G$ 
by giving a (set of) representative(s) in ${\rm GL}(2,F)$. 
\bigskip

 We write $T$ for the diagonal torus of $G$ and $B\supset T$ for the
 upper standard Borel subgroup. We denote by $T^0$ the maximal compact
 subgroup of $T$, i.e. the set of matrices with coefficients in
 $\ofr^{\times}$, and by $T^n$ the subgroup of matrices with
 coefficients in $1+\pfr^n$, $n>0$. 
\bigskip

 Let $k$, $l$ be integers satisfying $k+l\geqslant 0$. Then
 ${\mathfrak A}(k,l)= \matrice{\ofr}{\pfr^l}{\pfr^k}{\ofr}$ is an
 $\ofr$-order of ${\rm M}(2,F)$. We denote by $\Gamma_0 (k,l)$ the
 image in $G$ of its group of units. There are two conjugacy classes
 of maximal compact subgroups of $G$. The first one has representative
 $K=\Gamma_{0}(0,0)$. A representative $\tilde I$ of the second one
 is the semidirect  product of the  cyclic group generated by
 $\ds \Pi =\matrice{0}{1}{\pi}{0}$ with the  Iwahori subgroup 
$I=\Gamma_{0}(1,0)$. 
\bigskip

 The group $K$ is filtered by the normal compact open subgroups 
$$
K_n  =\matrice{1+\pfr^n}{\pfr^n}{\pfr^n}{1+\pfr^n} , \   n\geqslant 1
\ . 
$$
\noi The group $I$ is filtered by the normal compact subgroups $I_n$,
 $n\geqslant 1$, defined by:
$$
I_{2k+2}=\matrice{1+\pfr^{k+1}}{\pfr^{k+1}}{\pfr^{k+2}}{1+\pfr^{k+1}}\ ,
\ I_{2k+1}=\matrice{1+\pfr^{k+1}}{\pfr^{k}}{\pfr^{k+1}}{1+\pfr^{k+1}}\ ,
\ k\geqslant 0\ .
$$
\noi The subgroups $I_{n}$, $n\geqslant 1$, are normalized by
$\Pi$. 
\bigskip

 We denote by $X$ the Bruhat-Tits building of $G$. This is a uniform
 tree with valency $q+1$. As a $G$-set and as a simplicial complex $X$
  identifies  with the following complex. Its vertices are
 the homothety classes $[L]$ of full $\ofr$-lattices $L$ in the vector
 space $V=F^2$. Two vertices $[L]$ and $[M]$ define an edge is and
 only if there exists a basis $(e_1 ,e_2 )$ of $V$ such that, up to
 homothety, we have $L=\ofr e_1 \oplus \ofr e_2$ and $L=\ofr e_1
 \oplus \pfr e_2$.
\bigskip
 
The vertices of the standard apartment (i.e. the apartment stabilized
by $T$) are the $s_k =[\ofr\oplus\pfr^k ]$, $k\in \ZZ$. The element
$t_{\varpi}$ acts as $t_{\varpi} s_k =s_{k+1}$, $k\in \ZZ$. The maximal compact
subgroup $K$ is the stabilizer of $s_0$ and $\tilde I$ (resp. $I$) is the global
stabilizer (resp. pointwise stabilizer)  of the edge $[s_0 ,s_1
]$. If $l\geqslant k$, the pointwise stabilizer of the segment $[s_k
  ,s_l ]$ is $\Gamma_{0}(l,-k )$.

\section{Combinatorics of $\tdX_n$}

 We recall the construction of the directed graphs $\tdX_n$,
 $n\geqslant 1$.

For any integer $k\geqslant 1$, an {\it oriented $k$-path} in $X$ is
an injective sequence of vertices $(s_i )_{i=0,\dots ,k}$ in $X$ such
that, for $k=0,\dots, k-1$, $\{ s_i ,s_{i+1}\}$ is an edge in $X$. We
shall allow the index $i$ to run over any interval of integers of
length $k+1$. Let us fix an integer $n\geqslant 1$. The directed graph
 $\tdX_n$ is constructed as follows. Its edge set (resp. vertex set)
is the set of oriented $(n+1)$-paths (resp. $n$-paths) in $X$. If $a
=\{ s_0, s_1 ,...,s_{n+1}\}$ is an edge of $\tdX_n$, its head
(resp. its tail) is $a^+ =\{ s_1 ,s_2 ,...,s_{n+1}\}$ (resp. $a^- =\{
s_0, s_1 , ...,s_n\}$). The graphs we obtain this way are actually
simplicial complex. The group $G$ acts on $\tdX_n$ is an obvious way;
the action preserves the structure of directed graph.

When $n=2m$ is even, we have a natural simplicial projection 
$p=p_{n}~:$ $\tdX_{n}\rightarrow X$ given
on  vertices by $p(s_{-m},\dots ,s_0,\dots ,s_m )=s_0$. It is
$G$-equivariant. Let $e=\{ s_0 ,t_0\}$ be an edge of $X$. We are going
to describe the finite simplicial complex $p^{-1}(e)$. An edge in
$\tdX_n$ above the edge $e$ corresponds to an oriented $(2m+1)$-path of
one of the following forms:
\smallskip

i) $(s_{-m},s_{-m+1}, \dots ,s_{0},t_{0},\dots ,t_{m-1},t_{m})$

ii) $(t_{-m},t_{-m+1}, \dots ,t_{0},s_{0},\dots ,s_{m-1},s_{m})$
\smallskip

 Let $C_{2m-1}(e)$ the set of $(2m-1)$-paths $c=(s_{-m+1}, \dots ,s_0 ,t_0 ,\dots
 ,t_{m-1})$. We say that $c\in C_{2m-1}(e)$ {\it lies above} $e$.  Fix
 $c\in C_{2m-1}(e)$ and consider the simplicial sub-complex $\tdX_{2m}[e,c]$ of
 ${\tilde X}_{2m}$  whose edges correspond to the $(2m+1)$-paths of the form
$$
(a ,s_{-m+1}, \dots ,s_{0},t_{0},\dots ,t_{m-1},b)\, .
$$
\noi So $a$ (resp. $b$) can be any neighbour of $s_{-m+1}$
(resp. $t_{m-1}$) different from $s_{-m+2}$ (resp. $t_{m-1}$), with
the convention that $s_1 =t_0$ and $t_{-1}=s_0$. The simplicial
complex $\tdX_{2m}[e,c]$ is connected. It is indeed isomorphic to the
complete bipartite graph with sets of vertices:
$$
\{ a \ ; \ a\ {\rm neighbour\ of \ }s_{-m+1},\ a\not=s_{-m+2}\} \ {\rm and
} \ \{ b \ ; \ b\ {\rm neighbour\ of \ }t_{m-1},\ b\not=t_{m-2}\}\, .
$$

\smallskip

\noi {\bf Lemma 2.1}. {\it Let $e$ and $e'$ be two edges of $X$ and
$c\in C_{2m-1}(e)$, $c'\in C_{2m-1}(e')$. Then $\tdX_{e,c}\cap \tdX_{e',c'}\not=
\emptyset$ if and only if we are in one of the following cases:
\smallskip

i) $e=e'$ and $c=c'$ (so that $\tdX_{2m}[e,c]=\tdX_{2m}[e',c'])$;

ii) $e\cap e'$ is reduced to one vertex of $X$ and $c\cup c'$ is an
oriented $2m$-path in $X$. In that case $\tdX_{2m}[e,c]\cap \tdX_{2m}[e',c']$
is reduced to the vertex of $\tdX_{2m}$ corresponding to the $2m$-path
$c\cup c'$. }
\smallskip

\noi {\it Proof}. If $\tdX_{2m}[e,c]\cap \tdX_{2m}[e' ,c']\not=\emptyset$,
then $e\cap e' =p(\tdX_{2m}[e,c])\cap
p(\tdX_{2m}[e',c'])\not=\emptyset$. Assume first that $e=e'$. Then $c=c'$,
for if $c\not= c'$, then $\tdX_{2m}[e,c]\cap \tdX_{2m}[e',c']=\emptyset$;
indeed if $\tilde s$ is a vertex of $\tdX_{2m}[e,c]$ then it determines
$c$ uniquely. Now assume that $e\cap e'$ is a vertex. Let ${\tilde
s}\in \tdX_{2m}[e,c]\cap \tdX_{2m}[e',c']$. Then $\tilde s$ contains $c$ and
$c'$ as subsequences, with $c\not= c'$. So by a length argument
$s=c\cup c'$. Conversely if $c\cup c'$ is an oriented $2n$-path then
$c\cup c'$ is a vertex of $\tdX$ lying in $\tdX_{2m}[e,c]\cap     
\tdX_{2m}[e',c']$. 
\smallskip

\noi {\bf Corollary 2.2}. {\it For any edge $e$ of $X$, the connected
components of $p^{-1}(e)$ are the $\tdX_{2m}[e,c]$, where $c$ runs over
$C_{2m-1}(e)$.}

\medskip

Define a $1$-dimensional simplicial complex $Y_{2m-1}$ in the following
way. Its vertices are the connected components $\tdX_{2m}[e,c]$, where $e$
runs over the edges of $X$ and $c$ over $C_{2m-1}(e)$, and two vertices
$\tdX_{2m}[e,c]$ and $\tdX_{2m}[e',c']$ are linked by an edge if they
intersect. Note that $Y_{2m-1}$ is naturally a $G$-simplicial complex.  
\bigskip

\noi {\bf Corollary 2.3}. {\it As a $G$-simplicial complex $Y_{2m-1}$ is
canonically isomorphic to the complex $\tdX_{2m-1}$.}

\bigskip

 Assume that $m\geqslant 1$.  We say that an edge of
$\tdX_{2m-1}$ {\it lies above} a vertex $s_0$ of $X$ if as an 
oriented $2m$-path it
has the form $(s_{-m},\dots ,s_o,\dots ,s_{m})$. For any vertex $s_0$
of $X$ we write $\tdX_{2m-1}[s_o ]$ for the subsimplicial complex of $\tdY$
formed of the edges lying above $s_0$.  
\smallskip

\noi {\bf Lemma 2.5} {\it When $m=1$ the simplicial complexes
  $\tdX_{2m-1}[s_0 ] = \tdX_{1}[s_0 ] $
are connected.}
\smallskip

{\it Proof}. We may identify the neighbour vertices of $s_0$  in $X$ with the
points of the projective line ${\mathbb P}^{1}({\bar M})\simeq
{\mathbb P}^{1}({\mathbf k})$, where $s_{0}=[M]$ and ${\bar M}=M/\pk
M$. The vertices of $\tdX_{1}[s_0 ]$ are the oriented $1$-paths $(s_0
,x)$, $(y,s_0 )$, $x$, $y\in  {\mathbb P}^{1}({\bar M})$. Two
oriented $1$-paths of the form $(x,s_0 )$ and $(s_0 ,y)$ are linked by
the edge $(x,s_0,y)$. Let $(x,s_0 )$, $(y,s_0 )$ be two oriented
$1$-paths with $x\not= y$. Since $\vert {\mathbb P}^1 ({\mathbf
k})\vert \geqslant 3$, there exists $z\in {\mathbb P}^{1}({\bar M})$
distinct from $x$ and $y$. Then $(x,s_0 )$ is linked to $(s_o ,z)$ via
the path $(x,s_0 ,z)$ and $(s_0,z)$ is linked to $(y,s_0 )$ via the
path $(y,s_0 ,z)$. For vertices of the form $(s_0 ,x)$, $(s_0 ,y)$ the
proof is similar.
\bigskip

We now assume that $m>1$. We write $C_{2m-2}(s_0 )$ for the set
$(2m-2)$-paths of the form $(s_{-m+1},\dots ,s_0 ,\dots ,s_{m-1})$. For any
$c\in C_{2m-2}(s_0 )$, we consider the subsimplicial complex
$\tdX_{2n-1}[s_0 ,c ]$
of $\tdX_{2m-1}$ whose edges corresponds to the $2m$-paths of the form $(a,
s_{-n+1},\dots ,s_0 ,s_{n-1}, b)$. We have results similar to lemma
1.2, corollaries 2.2 and 2.3.
\smallskip

\noi {\bf Lemma 2.6}.  i) {\it For any vertex $s_0$ of $X$ and for
$c\in C_{2m-2}(s_0 )$,  $\tdX_{2m-1}[s_0 , c]$ is connected. It is indeed
isomorphic to a complete bipartite graph constructed on two sets of
$q=\vert {\mathbf k}\vert$ elements.}

\noi ii)  {\it Let $s$ and $s'$ be vertices of $X$, $c\in C_{2m-2}(s)$ 
and $c'\in C_{2m-2}(s')$. Then $\tdX_{2m-1}[s,c]\cap
\tdX_{2m-1}[s',c']\not= \emptyset$ if and only if
$s=s'$ and $c=c'$, or $\{s ,s'\}$ is an edge in $X$ and $c\cup c'$ is
an oriented $2n-1$-path. In this last case $\tdX_{2m-1}[s,c]\cap
\tdX_{2m-1}[s',c'] =\{ {\tilde s}\}$, where the vertex $\tilde s$ of
$\tdX_{2m-1}$ corresponds to the $(2n-1)$-path $c\cup c'$.

iii) For any vertex $s$ of $X$, the connected components of $\tdX_{2m-1}[s]$
are the $\tdX_{2m-1}[s,c]$, $c$ running over $C_{2m-2}(s)$. }
\bigskip

 We can consider the $1$-dimensional simplicial complex $Z_{2m-2}$ whose
 vertices are the connected components $\tdX_{2m-1}[s,c]$, $s$ running over
 the vertices of $X$ and $c$ over $C_{2m-2}(s)$, and where two connected
 components define an edge if and only if they intersect. Note that $Z_m$ is
 naturally a $G$-simplicial complex.
\smallskip

\noi {\bf Corollary 2.7}. {\it As a $G$-simplicial complex $Z_{2m-2}$
  is isomorphic to $X_{2n-2}$.}

\section{The cohomology of $\tdX_n$: first reductions}

If $\Sigma$ is a locally finite $1$-dimensional simplicial complex, we
write $\Sigma^0$ (resp. $\Sigma^{(1)}$,$\Sigma^1$) for its set of vertices (resp. non-oriented
edges, oriented edges). We let $C_0 (\Sigma)$ (resp. $C_1 (\Sigma)$) denote the
$\CC$-vector space with basis $\Sigma^0$ (resp. $\Sigma^1$). We define the space
 $C_{c}^{0}(\Sigma,\CC )=C_{c}^{0}(\Sigma)$ (resp.   $C_{c}^{1}(\Sigma,\CC
)=C_{c}^{1}(\Sigma)$) of oriented simplicial $0$-cochains
(resp. $1$-cochains) with compact support by:
\smallskip

$C_{c}^{0}(\Sigma)=$ space of all linear forms $f$~: $C_{0}(Z)\rightarrow
\CC$ such that $f(s)=0$ except for a finite number of vertices $s$;
\smallskip

$C_{c}^{1}(\Sigma)=$ space of all linear forms $\omega$~:
$C_{1}(\Sigma)\rightarrow \CC$ such that $f([a,b])=0$  except for a finite
number of oriented edges $[a,b]$ and $f([a,b])=-f([b,a])$. 
\smallskip

\noi We set $C_{c}^{k}(\Sigma)=0$ for $k\in \ZZ\backslash \{ 0,1\}$ and define a
coboundary map $d$~: $C_{c}^{0}(\Sigma)\rightarrow C_{c}^{1}(\Sigma)$ by
$df([a,b])=f(b)-f(a)$. The cohomology of the cochain complex
$(C_{c}^{\bullet}(\Sigma),d)$ computes the cohomology with compact support
$H_{c}^{\bullet}(\Sigma,\CC )=H_{c}^{\bullet}(\Sigma)$ of (the standard
geometric realization of) $\Sigma$. If $\Sigma$ is acted upon by a group $H$
whose action is simplicial then $(C_{c}^{\bullet}(\Sigma),d)$ is in a
straightforward way a complex of $H$-modules and its cohomology
computes $H_{c}^{1}(\Sigma)$ as a $H$-module.  When $T$ is finite we drop
the subscripts $c$.   
\smallskip

 Since the stablizer of a finite number of vertices of $X$ is open in
 $G$, we see that for $n\geqslant 1$,  the $G$-modules $C_{c}^{0}(\tdX_n )$,
 $C_{c}^{1}(\tdX_n )$ and therefore $H_{c}^{1}(\tdX_n )$ are
 smooth.

\medskip

 In the sequel we fix $m\geqslant 1$ and we abbreviate $\tdX_{2m}
 =\tdX$. The disjoint union $\displaystyle \tdX^1 = \bigsqcup_{e\in
 X^{(1)}} \tdX_{e}$, where $\tdX_{e}=p^{-1}(e)$, induces an
 isomorphism:
$$
\begin{array}{ccc}
C_{c}^{1}(\tdX )& \simeq & \bigoplus_{e\in X^{(1)}}C^{1}(\tdX_{e})\\
\omega & \mapsto & (\omega_{\vert C_{1}(\tdX_{e})})_{e\in X^{(1)}}
\end{array} \leqno (\hbox{3.1})
$$

\noi Similarly the non-disjoint union $\displaystyle
\tdX^{0}=\bigcup_{e\in X^{(1)}}\tdX_{e}^{0}$ induces an injection:
$$
\begin{array}{cccc}
j : & C_{c}^{0}(\tdX )& \hookrightarrow & \bigoplus_{e\in X^{(1)}}C^{0}(\tdX_{e})\\
 & f & \mapsto & (f_{\vert C_{0}(\tdX_{e})})_{e\in X^{(1)}}
\end{array} \leqno (\hbox{3.2})
$$

We have the following commutative diagram of $G$-modules:
$$
\begin{CD}
&      & H_{c}^{0}(\tdX )  @>>>  \bigoplus_{e\in
  X^{(1)}}H^{0}(\tdX_{e})  @>\varphi >>  {\rm coker j }     \\
 &   &   @VVV    @VVV    @|   \\
0  @>>> C_{c}^{0}(\tdX )  @>j>>  \bigoplus_{e\in
  X^{(1)}}C^{0}(\tdX_{e})  @>>>  {\rm coker j}  @>>>   0 \\
 &  &  @V{d}VV     @V{\oplus d_{e}}VV   @VVV  \\
 0  @>>>  C_{c}^{1}(\tdX )  @>>>  \bigoplus_{e\in X^{(1)}}
  C^{1}(\tdX_e )  @>>>  0  @>>>  0 \\
 & & @VVV @VVV @VVV \\ 
&      & H_{c}^{1}(\tdX )  @>>>  \bigoplus_{e\in
  X^{(1)}}H^{1}(\tdX_e )  @>>>  0   \\
\end{CD}
$$

\noi Here, for $e\in X^{(1)}$, $d_e$ denote the coboundary map
$C^{0}(\tdX_e )\rightarrow C^{1}(\tdX_e )$. Since $\tdX$ is connected
(\cite{[Br1]} Lemma 4.1) and non compact, we have $H_{c}^{0}(\tdX )=0$. So
{\it the snake lemma} gives the kernel-cokernel exact sequence:
$$
0\rightarrow \bigoplus_{e\in X^{(1)}}H^{0}(\tdX_{e}) \rightarrow {\rm
coker}j \rightarrow H_{c}^{1}(\tdX ) \rightarrow \bigoplus_{e\in
X^{(1)}}H^{1}(\tdX_e )\rightarrow 0
$$
\noi that is
$$
0\rightarrow {\rm coker j}/\varphi \big(\bigoplus_{e\in
X^{(1)}}H^{0}(\tdX_e )\big)\rightarrow H_{c}^{1}(\tdX ) \rightarrow \bigoplus_{e\in
X^{(1)}}H^{1}(\tdX_e )\rightarrow 0 \leqno (\hbox{3.3})
$$

Abbreviate $Y=Y_{2m-1}$.
\bigskip

\noi {\bf Lemma 3.4}. {\it We have a canonical isomorphism of
$G$-modules 
$${\rm coker j}/\varphi \big(\bigoplus_{e\in
X^{(1)}}H^{0}(\tdX_e )\big)\simeq H_{c}^{1}(Y )\, .
$$}
\smallskip

\noi {\bf Proof}. From corollary 2.2 we have 
$$
\bigoplus_{e\in X^{(1)}}C^{0}(\tdX_e )=\bigoplus_{e\in
X^{(1)}}\bigoplus_{c\in C(e)}C^{0}(\tdX_{e,c})\, .
$$
\noi So the map $j$ is given by $\displaystyle f\mapsto  \bigoplus_{e\in
X^{(1)}}\bigoplus_{c\in C(e)} f_{e,c}$, where $f_{e,c}=f_{\vert
C_{0}(\tdX_{e,c})}$. Consider the $G$-equivariant morphism of vector
spaces 
$$
\psi~:\ \bigoplus_{e\in X^{(1)}}\bigoplus_{c_\in
C(e)}C^{0}(\tdX_{e,c})\rightarrow C_{c}^{1}(Y )
$$
\noi given as follows. If $\sigma$ is an oriented edge of $Y$ then
their exist uniquely determined edges $e_o$, $e_{o}'$ of $X$, $c_o\in
C(e_o )$, $c_{o}'\in C(e_{o}')$, such that $\sigma$ corresponds to the
intersection $\tdX_{e_o ,c_o}\cap \tdX_{e_{o}',c_{o}'}=\{ s_o\}$,
$s_o\in \tdX^0$. We then set
$$
\psi [(f_{e,c})_{e,c}](\sigma )=f_{e_{o}',c_{o}'}(s_o )-f_{e_o ,c_o}(s_o )\, .
$$
\noi Then $\psi$ is surjective and its kernel is precisely
$j(C_{c}^{0}(\tdX ))$. So we may identify ${\rm coker}j$ with
$C_{c}^{1}(Y )$. From corollary 2.2, we have 
$$
\bigoplus_{e\in X^{(1)}}H^{0}(\tdX_e )=\bigoplus_{e\in
X^{(1)}}\bigoplus_{c\in C(e)}H^{0}(\tdX_{e,c})
$$
\noi so that we may identify $\bigoplus_{e\in X^{(1)}}H^{0}(\tdX_e )$
with $C_{c}^{0}(\tdY )$. Under our identifications the map $\varphi$
: $\bigoplus_{e\in X^{(1)}}H^{0}(\tdX_e ) \rightarrow {\rm coker}j$
corresponds to the coboundary map $d$~: $C_{c}^{0}(\tdY )\rightarrow
C_{c}^{1}(\tdY )$, and we are done since all our identifications are
$G$-equivariant. 
\smallskip

\noi {\bf Proposition 3.5}. {\it For $m\geqslant 1$, we have an
  isomorphism of $G$-modules: 
$$
H_{c}^{1}(\tdX_n )\simeq H_{c}^{1}(\tdX_{2m-1} ) 
\oplus \cind_{{\mathfrak K}_{e_o}}^{G}H^{1}(\tdX_{e_o})
$$
\noi for any edge $e_o$ of $x$ and where ${\mathfrak K}_{e_o}$ denotes
the stabilizer of $e_{o}$ in $G$. }
\smallskip

\noi {\it Proof}. From the short exact sequence (3.3) and lemma 3.4,
we have the exact sequence of $G$-modules:

$$
0\rightarrow H_{c}^{1}(Y )\rightarrow H_{c}^{1}(\tdX ) \rightarrow \oplus_{e\in
X^{(1)}}H^{1}(\tdX_e )\rightarrow 0 \leqno (\hbox{3.6})
$$

\noi Since $G$ acts transitively on the edges of $X$, $ \oplus_{e\in
X^{(1)}}H^{1}(\tdX_e )$ identifies with the compactly induced
representation $\cind_{{\mathfrak
K}_{e_o}}^{G}H^{1}(\tdX_{e_o})$. Moreover by [Vign ??]({\bf Trouver la
bonne r\'ef\'erence}) this  induced
representation is projective in the category of smooth complex
representations of $G$. So the sequence (3.7) splits. 

\medskip

We assume that $m\geqslant 1$ and we abbreviate $\tdX =\tdX_{2m-1}$. The disjoint
union $\displaystyle \tdX^{1}=\bigsqcup_{s\in X^0} \tdX_{s}^{1}$
induces an isomorphism:

$$
\begin{array}{ccc}
C_{c}^{1}(\tdX )& \simeq & \bigoplus_{s\in X^{0}}C^{1}(\tdX_{s})\\
\omega & \mapsto & (\omega_{\vert C_{1}(\tdX_{s})})_{s\in X^{0}}
\end{array} \leqno (\hbox{3.7})
$$

\noi Similarly the non-disjoint union $\displaystyle
\tdX^{0}=\bigcup_{s\in X^{0}}\tdX_{s}^{0}$ induces an injection:
$$
\begin{array}{cccc}
j : & C_{c}^{0}(\tdX )& \hookrightarrow & \bigoplus_{s\in
X^{0}}C^{0}(\tdX_{s})\\
 & f & \mapsto & (f_{\vert C_{0}(\tdX_{s})})_{s\in X^{0}}
\end{array} \leqno (\hbox{3.8})
$$

We have the following commutative diagram of $G$-modules:
$$
\begin{CD}
&      & H_{c}^{0}(\tdX )  @>>>  \bigoplus_{s\in
  X^{0}}H^{0}(\tdX_{s})  @>\varphi >>  {\rm coker j }     \\
 &   &   @VVV    @VVV    @|   \\
0  @>>> C_{c}^{0}(\tdX )  @>j>>  \bigoplus_{s\in
  X^{0}}C^{0}(\tdX_{s})  @>>>  {\rm coker j}  @>>>   0 \\
 &  &  @V{d}VV     @V{\oplus d_{s}}VV   @VVV  \\
 0  @>>>  C_{c}^{1}(\tdX )  @>>>  \bigoplus_{s\in X^{0}}
  C^{1}(\tdX_s )  @>>>  0  @>>>  0 \\
 & & @VVV @VVV @VVV \\ 
&      & H_{c}^{1}(\tdX )  @>>>  \bigoplus_{s\in
  X^{0}}H^{1}(\tdX_s )  @>>>  0   \\
\end{CD}
$$

\noi Here, for $s\in X^{0}$, $d_s$ denote the coboundary map
$C^{0}(\tdX_s )\rightarrow C^{1}(\tdX_s )$. By Lemma 2.4,  $\tdX$ is connected.
So  we have $H_{c}^{0}(\tdX )=0$ since $\tdX$ is non-compact.  The
{\it snake lemma} gives the kernel-cokernel exact sequence:

$$
0\rightarrow {\rm coker j}/\varphi \big(\bigoplus_{s\in
X^{0}}H^{0}(\tdX_s )\big)\rightarrow H_{c}^{1}(\tdX ) \rightarrow
\bigoplus_{s\in X^{0}}H^{1}(\tdX_s )\rightarrow 0 \leqno (\hbox{3.9})
$$

\noi {\bf Lemma 3.10}. {\it We have a canonical isomorphism of
$G$-modules 
$${\rm coker j}/\varphi \big(\bigoplus_{s\in
X^{0}}H^{0}(\tdX_s )\big)\simeq H_{c}^{1}(\tdX_{2m-2} )\, .
$$}
\smallskip

\noi {\it Proof}. It is similar to the proof of lemma 3.4 and relies
on lemma 2.6 and corollary 2.7.
\smallskip

\noi {\bf Proposition 3.11} {\it For $m\geqslant 1$, we have an isomorphism of
$G$-modules~:
$$
H_{c}^{1}(\tdX_{2m-1} )\simeq H_{c}^{1}(\tdX_{2m-2})\oplus \cind_{{\mathfrak
K}_{s_o}}^{G}H^{1}(\tdX_{s_o})
$$
\noi for any vertex $s_o$ and where ${\mathfrak K}_{s_o}$ denotes the
stabilizer of $s_o$ in $G$.}
\bigskip

\noi {\it Proof}. Similar to the proof of proposition 3.5.
%\smallskip

%\noi {\bf Proposition 2.12}. {\it Assume that $n>1$ and fix a vertex
%$s_o$ and an edge $e_o$ of $X$. We have an isomorphism of $G$-modules:
%$$
%H_{c}^{1}(\tdX_n )\simeq H_{c}^{1}(\tdX_{n-1})\oplus\cind_{{\mathfrak
%K}_{e_o}}^{G}H^{1}(\tdX_{e_o})\oplus \cind_{{\mathfrak
%K}_{s_o}}^{G}H^{1}(\tdX_{s_o}) \, .
%$$}

%\noi {\it Proof}. This is a consequence of proposition 2.5 and 2.11. 

\bigskip

 Recall \cite{[Br2]} that $\tX_0$ is different from $X$. This is a directed
 graph whose set of vertices is isomorphic to $X^0$ as a $G$-set and
 whose set of edges is isomorphic to the $G$-set of oriented edges of
 $X$.

\section{Determination of the inducing representations -- I}

 Let $m\geqslant 0$ be a fixed integer and $e_0 =[s_0 ,s_1 ]$ be the
 standard edge. The aim of this section is to determine the
 $\KK_{e_{0}}$-module $H^1 (\tdX_{2m} [e_0 ])$. Here we have
 $\KK_{e_0} ={\tilde I}$, the normalizer in $G$ of the standard Iwahori
 subgroup. We have the semidirect products:
$$
{\tilde I}=\langle \matrice{0}{1}{\varpi}{0}\rangle \ltimes I =E^{\times} I
$$
\noi for any totally ramified subfield extension $E/F\subset {\rm
  M}(2,F)$ such that $E^{\times}$ normalizes $I$.
\medskip

  We first assume that $m\geqslant 1$. By Corollary (2.2), we have 
the disjoint union:
$$
\tdX_{2m}[e_0 ]=\coprod_{c\in C_{2m-1}(e_0 )}\tdX_{2m} [e_0 ,c]\ .
$$
\noi The group $\tilde I$ acts transitively on $C_{2m-1}(e_0  )$. This
comes form the standard fact that $I$, the pointwise stabilizer of
$e_0$ acts transitively on the apartments of $X$ containing $e_0$. 

 Let $c_0\in C_{2m-1}(e_0 )$ be the path
$$
s_{-m+1}, ..., s_0 , s_1 ,...,s_m \ .
$$
\noi The global stabilizer of $\tdX_{2m}[e_0 ,c_0 ]$ in $\tilde I$ is
the pointwise stabilizer of $c_0$ in $\tilde I$, that is
$$
\Gamma_0
(m,m-1)=\matrice{\ofr^{\times}}{\pfr^{m-1}}{\pfr^m}{\ofr\times}=T^0
I_{2m-1}\ .
$$
\noi It follows that
$$
H^1 (\tdX_{2m}[e_0 ])={\rm ind}_{T^0 I_{2m-1}}^{\tilde I}H^1
(\tdX_{2m}[e_0 ,c_0 ])\ .
\leqno{(\hbox{4.1})}
$$

On the other hand, an easy calculation shows that the pointwise
stabilizer of $\tdX_{2m}[e_ 0,c_0 ]$ is $T^1 I_{2m}$, where $T^1$ is
the congruence subgroup of $T$ given by
$$
T^1 =\matrice{1+\pfr}{0}{0}{1+\pfr}\ .
$$
\noi So the $T^0 I_{2m-1}$-module $H^1 (\tdX_{2m}[e_0 , c_0 ]$ may be
viewed as a representation of the finite group $T^0 I_{2m-1}/T^1
I_{2m}$, that is a semidirect product of the cyclic group
$\bk^{\times}$ with the abelian group $I_{2m-1}/I_{2m}\simeq \bk
\oplus \bk$. 
\medskip

 Set $\Gamma =\tdX_{2m}[e_0 ,c_0 ]$. This is a finite directed
 graph. Let $\Sigma_{-m}$ (resp. $\Sigma_{m+1}$) denote the set
of verticed of $X$ that are neighbours of $s_{-m+1}$ and different from
$s_{-m+2}$ resp. neighbours of $s_m$ and different from
$s_{m-1}$. Then the vertex set of $\Gamma$ is
$$
\begin{array}{lll}
\Gamma^0 & =  & \{ (a,s_{-m+1},...,s_0 ,...,s_m )\ ; \ a\in
\Sigma_{-m}\}\coprod \{ (s_{-m+1}, ...,s_0,...,s_m, b)\ ; \ b\in
\Sigma_{m+1}\}\\
  &  \simeq  & \Sigma_{-m}\coprod \Sigma_{m+1}
\end{array}
$$

\noi and its edge set is

$$
\Gamma^1 =\{ (a,s_{-m+1}, ...,s_0,...,s_{m},b)\ ; \ a\in \Sigma_{-m},
\ b\in \Sigma_{m+1}\}\simeq \Sigma_{-m}\times \Sigma_{m+1}\ .
$$
\noi In particular $\Gamma$ is a bipartite graph based on two sets of
$q$ elements. In particular, its Euler character is given by
$$
\chi (\Gamma )=1-{\rm dim}_{\CC}H^{1}(\Gamma )=2q-q^2\ ,
$$
\noi so that 
$$
{\rm dim}_{\CC}H^1 (\Gamma )=q^2 -2q +1 =(q-1)^2\ . \leqno{(\hbox{4.2})}
$$

 Let $\CC [\Gamma^1 ]$ be the space of complex fonction on $\Gamma^1$
 and $\HH (\Gamma )$ be the space of harmonic $1$-cochains
on $\Gamma$:
$$
\HH (\Gamma )=\{ f\in \CC [\Gamma ]\ ; \ \sum_{a\in \Gamma^1 ,\ s\in
  a}[a:s]f(a)=0\ , \ \text{all } s\in \Gamma^0\}\ .
$$
\noi Here $[a:s]$ denote an incidence number. In our case :
$$
f\in \HH (\Gamma )\text{ iff }
\left\{
\begin{array}{lll}
\ds \sum_{a\in \Sigma_{-m}}f(a,s_{-m+1},...,s_m,b) & = & 0\text{, all }b\\
\ds \sum_{b\in \Sigma_{m+1}}f(a, s_{-m+1},...,s_m,b) & = & 0\text{, all }a
\end{array}
\right.
\leqno{(\hbox{Harm})}
$$

This is a standard result (see e.g. \cite{[Br2]}Lemma (1.3.2)), that, as
a $T^0I_{2m-1}/T^1 I_{2m}$-module, $H^1 (\Gamma )$ is isomorphic to
the contragredient module of $\HH (\Gamma )$. 

 An easy computation shows that we may identify $\Gamma^1$ with
 $\bk\times \bk$ in such a way that:
\medskip

\noi 1) an element of $\ds
I_{2m-1}=\matrice{1+\pfr^m}{\pfr^{m-1}}{\pfr^m}{1+\pfr^m}$ acts as
$$
\big( 1+\matrice{\varpi^m a}{\varpi^{m-1} b}{\varpi^{m}c}{\varpi^{m}d}
\big) .(x,y)=(x+\bar{b}, y+\bar{c})
$$
\noi for $a,b,c,d\in \ofr$, $x,y\in \bk$, and
\medskip

\noi 2) an element of $T^0$ acts as
$$
\matrice{a}{0}{0}{d}.(x,y)=(\bar{a}\bar{d}^{-1}x,\bar{d}\bar{a}^{-1}y)
$$
\noi and
\medskip

the condition (Harm) writes:
$$
f\in \HH (\Gamma ) \text{ iff }
\left\{
\begin{array}{lll}
\ds \sum_{x\in \bk}f(x,y) & = & 0, \text{ all } y\in \bk\\
\ds \sum_{y\in \bk}f(x,y) & = & 0, \text{ all } x\in \bk
\end{array}
\right.
$$
\noi A basis of $\CC [\Gamma ]$ is formed of the fonctions $\chi_1
\otimes \chi_2 (x,y)=\chi_1 (x)\chi_1 (y)$, where, for $i=1,2$,
$\chi_i$ runs over the characters of $(\bk ,+)$. It is clear that the
$(q-1)^2$ dimensional subspace of $\CC [\Gamma ]$ generated by the
$\chi_1 \otimes \chi_2$, $\chi_1 \not\equiv 1$, $\chi_2\not\equiv 1$,
is contained in $\HH (\Gamma )$. So using (4.2), we obtain:
$$
\HH (\Gamma )={\rm Span}\{ \chi_1 \otimes\chi_2\ ; \chi_i\in
\widehat{\bk^{\times}}, \chi_i\not\equiv 1\ ,
i=1,2\}\ . \leqno{(\hbox{4.3})}
$$

It follows from (4.3) that as an $I_{2m-1}/I_{2m}$-module, the space
$\HH (\Gamma )$ is the direct sum of $1$-dimensional representations
corresponding to the characters $\alpha = \alpha (\chi_1 ,\chi_2 )$,
$\chi_i \not\equiv 1$, $i=1,2$, given by
$$
\alpha \big(1+\matrice{\varpi^m a}{\varpi^{m-1}b}{\varpi^m
  c}{\varpi^{m}d}\big) = \chi_1 (b)\chi_2 (a)\ .
$$
\noi In particular $\HH (\Gamma )$ is isomorphic to its contragredient
and therefore isomorphic to $H^1 (\Gamma )$ as an
$I_{2m-1}/I_{2m}$-module.
In the language of strata (the reader may refer to \cite{[BH]}{\S}4), for
$\chi_i \not\equiv 1$, $i=1,2$, the character $\alpha (\chi_1 ,\chi_2
)$ corresponds to a stratum of the form $[{\mathcal I}, 2m,2m-1,\beta
]$, where $\mathcal I$ is the standard {\it Iwahori order} and
$\beta\in {\rm M}(2,F)$ is an element of the form $\ds
\Pi^{2m-1}\matrice{u}{0}{0}{v}$, $u,v\in \ofr^{\times}$.  
In the terminology of \cite{[BH]}{\S}4, page 98, this stratum is a {\it
  ramified simple stratum}. 

 We now have enough material to prove the following result.
\bigskip

\noi {\bf Proposition (4.4)}. {\it Let $\lambda$ be an irreducible
  constituent of 
$$
H^1 (\tdX_{2m}[e_0 ])={\rm  ind}_{T^0I_{2m-1}}^{\tilde I}
H^1 (\tdX_{2m}[e_0 ,c_0 ])\ .
$$
\noi Then the compactly induced representation ${\rm c-Ind}_{\tilde
    I}\lambda$ is irreducible, whence supercuspidal.}
\bigskip

\noi {\it Proof}. It is a standard result that an irreducible compactly induced
representation is supercuspidal (see \cite{[Mau]} or \cite{[Ca]}, page
194).

 The proof of the irreducibility is also standard by an argument due
 to Kutzko. But we repeat it for convenience. By Frobenius
 reciprocity,  the restriction of $\lambda$ to $I_{2m-1}$
 contains a character   $\alpha (\chi_1 ,\chi_2 )$ corresponding to a
 (ramified) simple stratum. Since $\lambda$ is irreducible and since
 $\tilde I$ normalizes $I_{2m-1}$, the restriction $\lambda_{\vert
   I_{2m-1}}$ is a direct sum $\alpha_1 \oplus \cdots \alpha_r$ of 
 $\tilde I$-conjugates of $\alpha (\chi_1 ,\chi_2
   )$. They all correspond to simple strata.  Let $g\in G$ be an
   element intertwining $\lambda$ with itself. Then by restriction it
   intertwines a character $\alpha_i$ with a character $\alpha_j$ for
   some $j=1,...,r$. By \cite{[BH]}{\S} Lemma (16.1), page 111, such an
   element $G$ must belongs to $\tilde I$. It follows that the
   $G$-intertwining of $\lambda$ is equal to $\tilde I$ and that the
   representation ${\rm c-Ind}_{\tilde I}^G \lambda$ is irreducible
   according to Mackey's irreducibility criterion (\cite{[Ca]}
   Proposition (1.5), page 195).
\bigskip

 We finally consider the case $m=0$. The directed graph $\tX_0$ has
 $X^0$ as vertex set. An edge $\{ t,s\}$ in $X$ gives rise to two
 edges $[s,t]$ and $[t,s]$ in $\tX_0$. Since the action of $G$ on
 $\tX_0$ preserves the structure of digraph, the $G$-module $H_{c}^1
 (\tX_0 )$ may be computed using the following complex :
$$
0\lra C_{c}^{0}(\tX_0 )\lra C_{c}^{(1)}(\tX_0 )
$$
\noi where $C_{c}^{(1)}(\tX_0 )$ is the space of (unoriented)
$1$-cochains, that is the space of maps from $\tX_0^{(1)}$ (unoriented
edges)  to $\CC$ with finite support. The coboundary map is here given
by $df[s,t]=f(t)-f(s)$. Consider the $G$-equivariant injection $j$~:
$C_{c}^{1}(X)\lra C_{c}^{(1)}(\tX_0 )$ given by $j(\omega)$~:
$[s,t]\mapsto \omega ([s,t])$. We have the commutative diagram of $G$-modules :

$$
\begin{CD}
0  @>>> C_{c}^{0}(X )  @>{\rm id}>>  C_{c}^{0}(\tX_0 )  @>>>  0  @>>>   0 \\
 &  &  @VVV     @VVV   @VVV  \\
 0  @>>>  C_{c}^{1}(X )  @>j>>  C^{1}(\tdX_0 )  @>>>
 C_{c}^{(1)}(\tdX_0 )/{\rm Im}j  @>>>  0 \\
\end{CD}
$$
\noi The quotient  $C_{c}^{(1)}(\tdX_0 )/{\rm Im}j$ identifies with
the subspace of $C_{c}^{(1)}(\tX_0 )$ formed of those functions $f$
satisfying $f([s,t])=f([t,s])$ for all edges $\{ s,t\}$ of $X$. This
subspace is nothing other than the compactly induced representation
${\rm c-Ind}_{\tilde I}^G {\mathbf 1}_{\tilde I}$. The cokernel exact
sequence writes:
$$
0\lra H_{c}^1 (X )\lra H_{c}^1 (\tX_0 ) \lra \cind_{\tilde I}^G
{\mathbf 1}_{\tilde I}\lra 0
$$
\noi Now we use the following two facts : 
\medskip

-- the representation ${\rm
  c-Ind}_{\tilde I}^G {\mathbf 1}_{\tilde I}$ is a projective object
of the category of smooth representations of $G$,
\smallskip

 -- the $G$-module $H_c^1 (X)$ is isomorphic to the Steinberg representation
 ${\mathbf S}{\mathbf t}_G$ of $G$ (\cite{[BoSe]})
\medskip

\noi to obtain:
\bigskip

\noi {\bf Proposition (4.5)}.{\it  The $G$-module $H_{c}^{1}( \tX_0 )$ is
isomorphic to ${\mathbf S}{\mathbf t}_G\oplus {\rm c-ind}_{\tilde I}^{G}{\mathbf
  1}_{\tilde I}$.}

\section{The inducing representations -- II}

 We now determine the
 $\KK_{s_0}$-module $H^1 (\tX_{2m+1}[s_0 ])$. The arguments are very
 often similar to those of the previous section and we will not give
 all details.  Since the case $m=0$
 requires slightly different techniques we postpone it to the end of the
 section and assume first  that $m>0$. 
\medskip

 Recall that the stabilizer $\KK_{s_0}$ of $s_0$ in $G$ is the image
 $K$ of ${\rm GL}(2, \ofr )$ in $G$. 
\medskip

  Let $c_0\in C_{2m}(s_0 )$ be the path $(s_{-m},...,s_0 ,...,s_m
  )$. Its pointwise stabilizer is $\Gamma_0 (m,m)=T^0 K_m$. So as a
  $K$-module, $H^1 (\tdX_{2m+1}[s_0 ])$ is isomorphic to the induced
  representation $\ds {\rm Ind}_{T^0 K_m}^K H^1 (\tdX_{2m+1}[s_0 ,c_0 ])$.
 Moreover the pointwise stabilizer of $\tX_{2m+1}[s_0 ,c_0 ]$ is $T^1
 K_{m+1}$ and $H^1 (\tdX [s_0 ,c_0 ])$ may be viewed as a
 representation of $T^0 K_m /T^1 K_{m+1}$. 
\medskip

 As in the previous section, one may consider the bipartite graph
 $\Omega$ whose both vertice sets identify with $\bk$, equiped with an
 action of $K_m$ on $\Omega^1$ given by
$$
\big[I_2 +\varpi^m \matrice{a}{b}{c}{d}\big] . (x,y)=(x+\bar{b}, y+\bar{c})\ , 
$$
\noi the action of $T^0$ being given by
$$
\matrice{a}{0}{0}{d}.(x,y)=(\bar{a}\bar{d}^{-1} x,
\bar{d}\bar{a}^{-1}y)\ .
$$
\noi Then the contragredient of the $T^0 K_m /T^1 K_{m+1}$-module $H^1 (\tdX [s_0 ,c_0 ])$ is
isomorphic to the space $\HH (\Omega )$ of harmonic cochains on
$\Omega$. As in the previous section this later space is generated by
the functions $\chi_1 \otimes \chi_2$, where $\chi_i$, $i=1,2$, runs
over the  non trivial characters of $(k,+)$. The line $\CC
\chi_1\otimes \chi_2$ is acted upon by $K_m$ via the character $\alpha
(\chi_1 ,\chi_2 )$ given by
$$
\alpha (\chi_1 ,\chi_2 )\big( I_2 +\varpi^m
\matrice{a}{b}{c}{d}\big) =\chi_1 (b)\chi_2 (a)\ .
$$
\noi It follows that $\HH (\Omega )$ is isomorphic to its
contragredient and that $H^1 (\tdX_{2m}[s_0 ,c_0 ])$ is the direct sum
of the characters $\alpha (\chi_1 ,\chi_2 )$, $\chi_i \not\equiv 1$,
$i=1,2$. 

 For $\chi_i \not\equiv 1$, $i=1,2$, the character $\alpha (\chi_1
 ,\chi_2 )$ corresponds to a stratum of  the form $[{\rm M}(2,\ofr ),
   m,m-1 ,\beta ]$, where $\beta\in {\rm M}(2,F)$ is given by 
$\ds \varpi^{-m}\matrice{0}{u}{v}{0}$, $u$, $v\in \ofr^{\times}$.
This stratum is either simple and non-scalar or split
fundamental
according to whether $uv$ mod $\pfr$  is a square in $\bk^{\times}$ or
not (here we have used the fact that ${\rm Char}\, (\bk )\not= 2$. 

 It is clear that $T^0$ leaves the set of characters corresponding to
 simple strata (resp. split fundamental strata)
 stable.
 So we may
 write
$$
H^1 (\tdX_{2m}[s_0 ,c_0 ])=H^1 (\tdX_{2m}[s_0 ,c_0 ])_{\rm
  simple}\oplus H^1 (\tdX_{2m}[s_0 ,c_0 ])_{\rm split}
$$
\noi where   $H^1 (\tdX_{2m}[s_0 ,c_0 ])_{\rm
  simple}$  (resp. $H^1 (\tdX_{2m}[s_0 ,c_0 ])_{\rm split}$) is the 
sub-$T^0 K_m$-module which decomposes as a $K_m /K_{m+1}$-module as a
direct sum of (characters corresponding to) simple non-scalar strata
(resp. split fundamental  strata). 

 We have a result similar to proposition (4.4), whose proof uses the
 same arguments.
\bigskip

\noi {\bf Proposition (5.1)}. {\it Let $\lambda$ be an irreducible
constituent of 
$$
{\rm Ind}_{T^0 K_m}^{K} H^1 (\tdX_{2m+1}[s_0 ,c_0 ])_{\rm simple}\subset
  H^{1}(\tdX_{2m+1}[s_0 ])\ .
$$
\noi Then the compactly induced representation $\cind_{K}^{G}\lambda$ is irreducible, whence supercuspidal.}
\bigskip

The study of ${\rm Ind}_{T^0 K_m}^{K} H^1 (\tdX_{2m+1}[s_0 ,c_0
])_{\rm split}$ is the aim of the next section.
\bigskip

 We are now going to determine the $K$-module structure of $H^1
 (\tdX_1 [s_0 ])$. Set $\GG = {\rm PGL}(2,\bk )\simeq K/K^1$ and write
 $\BB$ and $\TT$ for the upper Borel subgroup and diagonal torus of
 $\GG$ respectively. Let $\UU$ be the unipotent radical of $\BB$. 
 As a $K$-set the set of neighbour vertices 
of $s_0$ is isomorpic to ${\mathbb P}^1 (\bk )=\GG /\BB$.
\medskip
  
The graph $\Omega =\tdX_1 [s_0 ]$ has for vertex set the set of paths
of the form $(s,s_0 )$ or $(s_0 ,s)$ where $s$ runs over the neighbour
vertices of $s_0$ in $X$. So the space $C^0 (\Omega )$ of $0$-cochains
identifies with the space ${\mathcal F}({\mathbb P}^1 (\bk )\coprod
{\mathbb P}^1 (\bk ))$ of complex valued functions on the disjoint union
${\mathbb P}^1 (\bk )\coprod {\mathbb P}^1 (\bk)$. So has a
$\GG$-module $C^0 (\Omega )$ is isomorphic to ${\mathbf 1}_\GG \oplus
\St_\GG \oplus {\mathbf 1}_\GG \oplus \St_\GG$, where ${\mathbf 1}$
denotes a trivial representation and $\St$ a Steinberg representation.

\medskip

 The $\GG$-set $\Omega^1$ is the set of paths of the form $(s,s_0,t)$,
 where $s$ and $t$ are two different neighbour vertices of $s_0$. This
 $\GG$-set is isomorphic to the quotient $\GG /\TT$.  The space
 $C^{(1)}(\Omega )$ of unoriented $1$-cochains identifies as
 $G$-module with the space ${\mathcal F}(\GG /\TT)$.

Fix a
 non-trivial character $\psi$ of $\UU$. It is well knows that the
 induced representation ${\rm Ind}_\UU^\GG \psi$ is multiplicity
 free. Its irreducible constituent form by definition the {\it
   generic} (irreducible) representations of $\GG$. Moreover an irreducible
 representation is generic if ans only if it is not a character.

 We have a natural $G$-equivariant map $\Phi$~: ${\mathcal F}(\GG /\TT )\lra {\rm Ind}_\UU^\GG
 \, \psi$, given by 
$$
\Phi (f)(g)=\sum_{u\in \UU} f(gu){\bar \psi}(u)\ , \ f\in {\mathcal F}(\GG /\TT ),
\ g\in \GG \ .
$$
\noi If a  function $f$ lies in the kernel of $\Phi$, then we have
$\ds \sum_{u\in \UU}f(gu)\theta (u)=0$, for all $g\in \GG$ an all
non-trivial character $\theta$ of $\UU$. Indeed it suffices to use the
fact that the action of $\TT$ on $\UU$ by conjugation acts
transitively on the non-trivial characters of $U$ and the right
invariance of $f$ under the action of $T$. So the kernel of $\Phi$
consists of the function $f$ such that $u\mapsto f(gu)$ is constant
function on $U$, for all $g\in G$. In other words ${\rm Ker}\, \Phi =
{\mathcal F}(G/B)\simeq {\mathbf 1}_\GG\oplus \St_\GG$. By a dimension
argument, we see that $\Phi$ is surjective. It follows that 
$$
C^{(1)}(\Omega )\simeq {\rm Ind}_\UU^\GG \, \psi \oplus {\mathbf 1}_\GG
\oplus \St_\GG \ .
$$

We have the cochain complex of $G$-modules:
$$
0 \lra C^0 (\Omega )\lra C^{(1)}(\Omega )\lra 0
$$
\noi Since $\Omega$ is connected the kernel of the coboundary operator
is the trivial module $\CC$. Hence in the Grothendieck groups of
$G$-modules, we have: $dC^0 (\Omega )\simeq 2.{\mathbf 1}_\GG +2.\St_\GG
-{\mathbf 1}_\GG ={\mathbf 1}_\GG +2.\St_\GG$. Therefore
$$
H^1 (\Omega )=C^1 (\Omega )/dC^0 (\Omega )\simeq {\rm Ind}_\UU^\GG \psi
+{\mathbf 1}_\GG +\St_\GG -{\mathbf 1}_\GG -2.\St_\GG = {\rm
  Ind}_\UU^\GG \psi -\St_\GG\ .
$$

Since $q=\vert\bk\vert$ is odd, there exists a  unique
non-trivial character of $\bk^{\times}/(\bk^{\times})^2$, that we denote by $\chi_0$. 
The  irreducible constituents  
 of the Gelfand-Graev representation ${\rm
  Ind}_\UU^\GG \psi$ are the following:
\bigskip

 -- the irreducible cuspidal representations of $\GG$,
\smallskip

 -- the principal series ${\rm Ind}_{\BB}^{\GG}\chi\otimes \chi^{-1}$,
 where $\chi$~: $\bk^{\times}\lra \CC^{\times}$ is a character such
 that $\chi^2\not\equiv 1$ (i.e. $\chi\not\in \{ {\mathbf 1}, \chi_0 \}$). 
\smallskip

 -- the steinberg representation $\St_\GG$,
\smallskip

 -- (when $q$ is odd) the twisted representation $\St_\GG\otimes
 \chi_0$.
\bigskip

 If $\sigma$ is a cuspidal representation of $\GG = K/K^1$, then the
 induced representation $\cind_K^G \sigma$ is irreducible and
 supercuspidal (\cite{[BH]}, (11.5), page 81). Such a 
 representation of $G$ is called a {\it level $0$ supercuspidal} representation.
\bigskip

 A principal series of $\GG =K/K^1$ may be writen as ${\rm
   Ind}_{I}^{K}\rho$, where $\rho$ is a character of $I/I^1$. The pair
 $(I,\rho )$ is actually a type in the sense of Bushnell and Kutzko's type
 theory. For technical reason we postpone definitions and references
 to the next section. Since the representation ${\rm
   Ind}_{I}^{K}\rho$ is irreducible, it is a type for the same
 constituent as $(I,\rho )$. 
\bigskip

To sum up, we have proved the following.
\bigskip

\noi {\bf Proposition (5.2)}. {\it An irreducible constituent $\lambda$
  of $H^1 (\tdX_1 [s_0 ])$ is of one of the following forms}
\medskip

 (i)  {\it the inflation of a cuspidal representation of $\GG$; in that case
 $\cind_K^G\lambda$ is a level $0$ irreducible supercuspidal
 representation of $G$.}
\smallskip

 (ii) {\it  the inflation to $K$ of the representation
 $\St_\GG\otimes \chi_0$,}
\smallskip

 (iii) {\it  a type of the form ${\rm Ind}_{I}^{K}\rho$, where the $\rho$ is
 inflated from a character of $I/I^1\simeq (\bk^\times
 \times\bk^\times )/\bk^\times$ of the form $\chi\otimes \chi^{-1}$,
 $\chi^2\not\equiv 1$. }
\bigskip 

\noi Note that in (iii), the pair $(K, {\rm Ind}_{I}^{K}\, \rho )$ is
a principal series type.

\section{The inducing representations -- III}

 We keep the notation as in the previous section. To determine the
 structure of ${\rm Ind}_{T^0 K_m}^{K} H^1 (\tdX_{2m+1}[s_0
   ,c_0])_{\rm split}$,  we first recall crucial facts on split strata and
 types for principal series representations. The basic reference for
 {\it type theory} is \cite{[BKtypes]}.
\bigskip

 Let $\chi$ be a character of $T$, that we view as a character of
 $T^0$ by restriction. Assume that the conductor of $\chi$ is $n>0$ : 
$T^n \subset {\rm Ker}\, \chi$ and $n$ is minimal for this
 property. Set
$$
J_{\chi}=\matrice{\ofr^{\times}}{\ofr}{\pfr^n}{\ofr^{\times}}
=\Gamma_{0}(\pfr^n )\ .
$$
\noi If $U$ and $\bar U$ denotes the groups of upper and lower
unipotent matrices respectively, then $J_{\chi}$ has an Iwahori
decomposition:
$$
J_{\chi}= (J_{\chi}\cap {\bar U}).(J_{\chi}\cap T ).(J_{\chi}\cap U)
$$
\noi and one may define a character $\rho_{\chi}$ of $J_{\chi}$ by
$$
\rho_{\chi}(\bar{u}t^0 u)=\chi (t^0 )\ , \ \bar{u}\in  J_{\chi}\cap
    {\bar U} , \ u\in J_{\chi}\cap U , \ t^0 \in T^0\ .
$$

Let ${\mathcal R}_{[T,\chi ]}$ be the Bernstein component of the
category of smooth representations  of $G$ whose objects are the representations
$\mathcal V$ satisfying the following property : any irreducible
subquotient of $\mathcal V$ occurs in a
 parabolically  induced representation ${\rm
  Ind}_{\mathcal B}^G (\chi\otimes \chi_0 )$, where $\mathcal B$ is
 a Borel subgroup with Levi component $T$ and $\chi^0$ an unramified
 character of $T$. We then have.
\bigskip

\noi {\bf Theorem (6.1)} (A. Roche) {\it The pair $(J_{\chi},
  \rho_{\chi})$ is a type for ${\mathcal R}_{[T,\chi ]}$.}
\bigskip

This is indeed Theorem (7.7) of \cite{[Ro]}. Note that our $J_{\chi}$
is not exactly Roche's one, but a conjugate under an element of $T$
(see \cite{[Ro]}, Example (3.5)). 
\bigskip

\noi {\bf Proposition (6.2)}. {\it With the notation as before, assume
  that $\chi_{\vert T^0}$ is not of the form $\alpha \circ {\rm Det}$,
  where $\alpha$ is a character of $\ofr^{\times}$ (necessarily of
  order $2$). Then the  induced representation $\ds {\rm
    Ind}_{J_{\chi}}^{K}\rho_{\chi}$ is  irreducible. In particular
 it is a type for ${\mathcal R}_{[T,\chi ]}$.}
\bigskip

\noi {\it Proof}. Let $W$ be the extended affine Weyl group of $G$
w.r.t. $T$ and set $W_{\chi}=\{ w\in W\ ; \ w\chi =\chi\}$. Then by
Theorem (4.14) of \cite{[Ro]}, the $G$-intertwining of $\rho_{\chi}$
is $J_{\chi}W_{\chi}J_{\chi}$. The hypothesis on $\chi$ forces
$W_{\chi} = T/T^0$. So $(J_{\chi}W_{\chi}J_{\chi})\cap K =J_{\chi} T^0
J_{\chi}=J_{\chi}$, and we may apply Mackey's criterion of irreducibility.
 
\bigskip

 For $n>0$ and $q\in \{ 0,...,n\}$, define compact open subgroups of
 $G$ as follows:
$$
{}_q\hh_1 =\matrice{1+\pfr^n}{\pfr^q}{\pfr^{n+1}}{1+\pfr^n}\text{ and }
  {}_q\hh_2 =\matrice{1+\pfr^{n+1}}{\pfr^q}{\pfr^{n+1}}{1+\pfr^{n+1}}
\ .
$$
\noi These groups are particular cases of groups considered in
\cite{[BrCrelle]}, {\S}(2.3). The quotients ${}_q\hh_1 /{}_q\hh_2$,
$q=0,...,n$, are abelian, and for $\alpha\in \bk^{\times}$, one may define a
character $\psi_{\alpha}$ of ${}_q\hh_1 /{}_q\hh_2$ by the formula:
$$
\psi_{\alpha}\big( I_2+\matrice{\varpi^{n}a}{\varpi^q
  b}{\varpi^{n+1}c}{\varpi^{n}d}\big) ={\mathbf \psi}(\alpha (a-d))
$$
\noi where $\mathbf \psi$ is a fixed non-trivial character of $(\bk ,+)$.
In fact, $(\psi_{\alpha})_{\vert {}_n\hh_{1}}$ is  the
restriction to  ${}_n\hh_1$ of a split 
fundamental stratum of $K_n /K_{n+1}$. We shall need the following result.
\bigskip

\noi {\bf Lemma (6.3)}. {\it  If a smooth representation of $K$
  contains $(\psi_{\alpha})_{\vert {}_n\hh_1}$ by restriction, then it
  contains the character $(\psi_{\alpha})_{\vert {}_0\hh_1}$.}
\bigskip

\noi {\it Proof}. Since the characteristic of $\bk$ is not $2$, then
$\alpha\not= -\alpha$ 
$(\psi_{\alpha})_{\vert {}_n\hh_1}$ is the restriction to ${}_n\hh_1$
of a split fundamental stratum of $K_n /K_{n+1}$. Our lemma is then a
particular case of \cite{[BrCrelle]}, Lemma (2.4.5).
\bigskip

% We now split the $T^0 K_m$-module $H^{1}(\tdX_{2m+1}[s_0 ,c_0
% ])_{split/scalar}$ as
%$$
%H^{1}(\tdX_{2m+1}[s_0 ,c_0 ])_{split}\oplus  
%H^{1}(\tdX_{2m+1}[s_0 ,c_0 ])_{scalar}
%$$
%\noi where $H^{1}(\tdX_{2m+1}[s_0 ,c_0 ])_{split}$ is the submodule
%where $K_m /K_{m+1}$ acts via split strata and 
%$H^{1}(\tdX_{2m+1}[s_0 ,c_0 ])_{scalar}$ the submodule where $K_m
%/K_{m+1}$ acts via essentially scalar strata. The first module is zero
%when the characteristic of $\bk$ is $2$, while the second is zero when
%the characteristic of $\bk$ is not $2$.
%\bigskip

\noi {\bf Proposition (6.4)}. {\it  Let $\lambda$ be an irreducible
  constituent of $\ds {\rm Ind}_{T^0 K_m}^K H^{1}(\tdX_{2m+1}[s_0 ,c_0
  ])_{split}$. Then with the notation as above,  $\lambda$ is of the
  form  ${\rm Ind}_{J_{\chi}}^{K}\rho_{\chi}$, for some principal
  series type $(J_{\chi},\rho_{\chi})$ with $\chi$ of conductor
  $m+1$.}
\bigskip

\noi {\it Proof}. We know that such a $\lambda$ contains a split
fundamental stratum of the form $[{\rm M}(2,\ofr ), m, m-1 ,b]$, where
$\ds b=\varpi^{-m}\matrice{0}{u}{v}{0}$, $u,v\in \ofr^{\times}$, and
$uv$ is a square {\it modulo} $\pfr$. If $\alpha\in \ofr$ is such that
$\alpha^2 \equiv uv$ {\it mod} $\pfr$, then the stratum is equivalent
to a $K$-conjugate of $[{\rm M}(2,\ofr ),m,m-1 ,b']$, where $\ds 
b'=\varpi^{-m}\matrice{\alpha}{0}{0}{-\alpha}$. So we deduce that
$\lambda$ contains this latter stratum by restriction. Now consider
the group ${}_q\hh_1$ for $n=m$. The representation $\lambda$ contains
the character $(\psi_{\alpha})_{\vert {}_{n}\hh_1}$ by restriction. By
applying Lemma (6.3) we obtain that it contains the character
$(\psi_{\alpha})_{\vert {}_0\hh_1 }$. This character  clearly
  extends to  $T^0{}_0\hh_1 =\Gamma_{0} (m+1, 0)$ and the quotient
  $T^{0}{}_0\hh_1 /{}_0\hh_1$ is abelian. It follows that $\lambda$
  contains and extension of $\psi_{\alpha}$ to $\Gamma_0
  (m+1,0)$. Such an extension is of the form $(J_{\chi},\rho_{\chi})$,
  for some character $\chi$ of $T$ of conductor $m+1$. The fact that
  $\lambda$ is induced from $(J_{\chi}, \rho_{\chi})$ follows from
  Proposition (6.2).

\section{Synthesis}

We now prove Theorems A and B of the introduction.
\bigskip

 By Proposition (3.5) and (3.11), we have isomorphisms of $G$ modules :

\begin{equation}
H_c^1 (\tdX_{2m})\simeq H_c^1 (\tdX_{2m-1}) \oplus \cind_{\KK_1}^G \, H^1 (\Sigma_{2m}), \ m\geqslant 1 .
\end{equation}

\begin{equation}
H_c^1 (\tdX_{2m=1})\simeq H_c^1 (\tdX_{2m}) \oplus \cind_{\KK_0}^G \, H^1 (\Sigma_{2m+1}), \ m\geqslant 0 .
\end{equation}

 Recall that with the notation of the introduction, we have :
\medskip

 -- $\Sigma_{2m}=(\tdX_{2m})_{e_0}$, $\Sigma_{2m+1}=(\tdX_{2m+1})_{s_0}$,
\smallskip

 -- $\KK_0 ={\mathcal K}_{s_0}$, $\KK_1 ={\mathcal K}_{e_0}$.
\bigskip

 Moreover, by Proposition (4.5), we have

\begin{equation}
H_c^1 (\tdX_0 )\simeq {\rm St}_G \oplus \cind_{\KK_1}^G \, H^1 (\Sigma_0 )
\end{equation}

\noi so that (1) holds for $m=0$. Hence Theorem A follows from (1) and 
(2) by a straightforward inductive argument. 
\bigskip

Theorem B follows from the discription of the irreducible components
 of $H^1 (\Sigma_n )$ given in Proposition (4.4) ($n$ even and $n>0$), 
Proposition (4.5) ($n=0$), and Propositions (5.1) and (6.4) ($n$ odd).

%%%%%%%%%%%%%%%%%%%%%%%%%%%%%%%%%%%%%%%%%%%%%%%%%%%%%%%%%%%%
%%%  Bibliographie %%%
%%%%%%%%%%%%%%%%%%%%%%

\bigskip

\bigskip

\begin{center}
Laboratoire de Math\'ematiques et

UMR 7348 CNRS
\medskip

SP2MI - T\'el\'eport 2

Bd M. et P. Curie BP 30179

86962 Futuroscope Chasseneuil Cedex

France
\medskip

E-mail~: paul.broussous@math.univ-poitiers.fr
\end{center}


\begin{thebibliography}{99}

\bibitem{[BrCrelle]} P. Broussous, {\it Minimal strata for ${\rm
    GL}(m,D)$}, J. reine angew. Math. {\bf 514} (1999), 199-236.


\bibitem{[Br1]} P. Broussous, {\it Simplicial complexes lying
  equivariantly over the affine building of ${\rm GL}(N)$},
  Math. Annalen, {\bf 329} (2004), 495--511.

\bibitem{[Br2]} P. Broussous, {Representations of ${\rm PGL}(2)$ of a
  local field and harmonic cochains on graphs}, 
  Ann. Fac. Sci. Toulouse Math. (6)  18  (2009),  no. 3, 495--513.

%\bibitem{[Brown]} K.S. Brown, {\it Building}, Springer Verlag, 1989.

\bibitem{[BH]} C.J. Bushnell and G. Henniart, {\it The local Langlands
  conjecture for ${\rm GL}(2)$}, Grundlehren des Math. Wiss. vol. 335,
  Springer, (2006).

\bibitem{[BKtypes]} C.J. Bushnell and P.C. Kutzko, {\it Smooth
  representations of reductive $p$-adic groups; structure theory via
  types}, Proc. London Math. Soc. (3)  77  (1998),  no. 3, 582--634.  


%\bibitem{[BT]} F. Bruhat and J. Tits, {\it Groupes r\'eductifs sur un
%corps local. I. Donn\'ees radicielles valu\'ees}, Publ. Math. IHES
%{\bf 41} 1972, 5-251.

\bibitem{[BK]} C.J. Bushnell and P.C. Kutzko, {\it Smooth
representations of reductive $p$-adic groups: structure theory via
types}, Proc. London Math. Soc. (3) {\bf 77} 1998, 582-634.


\bibitem{[BoSe]} A.  Borel and J.--P. Serre, {\it  Cohomologie \`a supports
 compacts des immeubles de Bruhat-Tits; applications à la cohomologie 
des groupes $S$-arithm\'etiques},  C. R. Acad. Sci. Paris S\'er. A-B 272 1971 A110–A113. 


\bibitem{[Ca]} H. Carayol, {\it Repr\'esentations cuspidales du groupe
lin\'eaire}, Ann. Sci. Ecole Norm. Sup. (4) {\bf 17} no. 2 1984,
191-225.

%\bibitem{[Cass1]} W. Casselman, {\it On some result of Atkin and
%Lehner}, Math. Ann. {\bf 201}, 1873, 301-314.

%\bibitem{[Cass2]} W. Casselman, {\it The restriction of a
%representation of ${\rm GL}_2 (k)$ to ${\rm GL}_2 ({\mathfrak o})$},
%Math. Ann. {\bf 206}, 1973, 311-318.

%\bibitem{[JPSS]} H. Jacquet, I. Piateski-Shapiro, J. Shalika, {\it
%Conducteur des repr\'esentations du groupe lin\'eaire},
%Math. Ann. {\bf 256} no. 2, 1981, 199-214.

%\bibitem{[He]} G. Henniart, {\it Sur l'unicit\'e des types pour ${\rm
%    GL}_2$}, appendix to {\it Multiplicit\'e modulaires et
%    repr\'esentations de ${\rm GL}(2,{\mathbb Z}_p )$ et ${\rm Gal}({\bar
%      {\mathbb Q}}_p/{\mathbb Q}_p)$ en $l=p$}, by C. Breuil and A. M\'ezard.

\bibitem{[Ku]} P.C. Kutzko, {\it On the supercuspidal representations
of ${\rm GL}(2)$, I, II}, Amer. J. Math. Vol. {\bf 100}, 1978, 43-60
and 705-716.

\bibitem{[Mau]} F.I. Mautner, {\it Spherical functions over $\mathfrak
 P$-adic fields. II}, Amer. J. Math. {\bf 86} (1964), 171-200.

\bibitem{[Ro]} A. Roche, {\it Types and Hecke algebras for principal
series of split reductive $p$-adic groups}, Ann Sci. Ecole
Norm. Sup. (4) {\bf 31} no. 3 1998 361-413. 

%\bibitem{[SS]} P. Schneider and U. Stuhler, {\it Representation theory
%and sheaves on the Bruhat-Tits building}, Publ. Math. IHES {\bf 85}
%1997, 97-191.

%\bibitem{[Wg]} J.B. Wagoner, {\it Homotopy theory for the $p$-adic
%special lineal group}, Comment. Math. Helvetici {\bf 50} 1975, 535-559.

\end{thebibliography}
\end{document}